\documentclass{elsart}
\usepackage{array,amsmath,amstext,amssymb,amsbsy,graphicx}
\usepackage{pdfsync}
 \usepackage{color}
\usepackage{txfonts}

\newcommand{\bfG}{{\boldsymbol G}}
\newcommand{\bfT}{{\boldsymbol T}}
\newcommand{\Gr}{\text{Grad}}

\newcommand{\bfxi}{\boldsymbol \xi}
\newcommand{\bfn}{\boldsymbol n}
\newcommand{\bfN}{\boldsymbol N}
\newcommand{\bfB}{\boldsymbol B}
\newcommand{\bfkappa}{\boldsymbol\varkappa}
\newcommand{\bfF}{\boldsymbol F}
\newcommand{\bfL}{\boldsymbol L}

\newcommand{\bff}{\boldsymbol f}

\newcommand{\bfJ}{\boldsymbol J}
\newcommand{\bfu}{\boldsymbol u}

\newcommand{\bfX}{\boldsymbol X}

\newcommand{\bfeps}{\boldsymbol\varepsilon}
\newcommand{\bfsig}{\boldsymbol\sigma}
\newcommand{\bfv}{\boldsymbol v}
\newcommand{\bfw}{\boldsymbol w}

\newcommand{\bfI}{\boldsymbol I}
\newcommand{\bfP}{\boldsymbol P}
\newcommand{\bfx}{\boldsymbol x}

\newcommand{\bfzero}{\boldsymbol 0}

\begin{document}
\begin{frontmatter}
\title{Tangential differential calculus and the finite element modeling of a large deformation elastic membrane shell problem}
\author[PH]{Peter Hansbo,} 
\author[MGL]{Mats G. Larson,}
\author[FL]{Fredrik Larsson}
\address[PH]{Department of Mechanical Engineering, J\"onk\"oping University,
SE-55111 J\"onk\"oping, Sweden}
\address[MGL]{Department of Mathematics and Mathematical Statistics,
Ume{\aa} University, SE--901 87 Ume{\aa}, Sweden}
\address[FL]{Department of Applied Mechanics, Chalmers University of Technology, SE--412 96 G\"{o}teborg, Sweden}
\date{}

\maketitle
\begin{abstract}
We develop a finite element method for a large deformation membrane elasticity problem on
meshed surfaces
using a tangential differential calculus approach that avoids the use
of classical differential geometric methods.
The method is also applied to form finding problems.
\end{abstract}
\end{frontmatter}

\section{Introduction}

In \cite{HaLa14}, we introduced a finite element method for the solution of the linear membrane shell problem.
This method was based on the use of tangential differential calculus which allows the use of Cartesian coordinates
for establishing the discrete system and for representing the displacements. The shell model underpinning the finite element method
is equivalent to the classical one given, e.g., by Ciarlet \cite{Ci00}, as demonstrated by Delfour and Zol\'esio \cite{DeZo97}.
We showed in \cite{HaLa14} that our finite element method can be viewed as a generalization of the classical flat facet element method (discussed, e.g., in Chapelle and Bathe \cite{ChBa03}), allowing
for higher order approximations on curved elements. Similar methods were also used in the context of beams in \cite{HaLaLa14}.

In this paper we extend the method of \cite{HaLa14} to the case of large deformation hyperelasticity. Our method is Lagrangian and uses the first Piola--Kirchhoff stress tensor to represent the stress field. Applications are given to
a general membrane problem using the Mooney--Rivlin model as an example of a constitutive law, and also
to \emph{form finding}\/, designed to determine structural shape
of compressed lightweight structures by inversion of tensile structures in the form of hanging models, cf., e.g., \cite{BlWuDaCa05}.
We focus on area minimization, a basic form finding method.

An approach similar to the one we propose was used in the context of interfaces by Monteiro, He, and Yvonnet \cite{MoHeYv11},
with the difference that the material models in \cite{MoHeYv11} were developed directly in plane stress, whereas we
use three dimensional modeling. 

The paper is organized as follows: in Section \ref{largedef} we introduce our formalism and discuss the large deformation problem we wish to solve; in Section \ref{FEM} we introduce and discuss the finite element method and our iterative solution method. In Section \ref{constitutive} we briefly discuss our examples of constitutive laws, and in Section \ref{numerical} we give some numerical examples.

\section{The large deformation membrane shell problem\label{largedef}}

\subsection{Basic notation}

We begin by recalling the fundamentals of the tangential calculus approach following Delfour and Zolesio \cite{DeZo95,DeZo97}. Let the shell in its undeformed configuration, in a Cartesian coordinate system $\bfX$, occupy a smooth two-dimensional surface $\Gamma$ embedded in ${\mathbb{R}}^3$, with outward pointing normal $\bfN$.
We denote the signed distance function relative to $\Gamma$ by $D(\bfX)$, for $\bfX\in \Bbb{R}^3$, fulfilling
\begin{equation}\label{Normal}
\left.\frac{\partial D(\bfX)}{\partial\bfX}\right|_\Gamma = \bfN(\bfX) .
\end{equation}
We can define the initial, undeformed, domain occupied by the membrane by
\[
\Omega_t = \{\bfX\in \Bbb{R}^3: \vert D(\bfX) \vert < t/2\},
\]
where $t$ is the thickness of the membrane.
The closest point projection $\bfX_\Gamma(\bfX):\Omega_t \rightarrow \Gamma$
is given by
\begin{equation}\label{closestP}
\bfX_\Gamma = \bfX -D(\bfX)\nabla D(\bfX), \quad\nabla D:= \frac{\partial D}{\partial \bfX} ,
\end{equation}
and thus the linear projector
$\bfT_\Gamma = \bfT_\Gamma(\bfX_\Gamma)$, onto the tangent plane of $\Gamma$, is given by
\[
\bfT_\Gamma := \bfI -\bfN\otimes\bfN ,
\]
where $\otimes$ denotes the outer product, and we can define the surface gradient
$\nabla_\Gamma$ as
\begin{equation}
\nabla_\Gamma := \bfT_\Gamma\cdot\nabla  .
\end{equation}
Note in particular that $\bfT_\Gamma$ is idempotent.

\subsection{Surface displacements}

We assume that the placement $\bfx$ of the original domain $\Omega_t$ can be given in terms of the original configuration $\bfX$ as
\begin{equation}\label{kinematics}
\bfx = \bfX +\bfu (\bfX_\Gamma) + \bar{\bfu}(\bfX_\Gamma)D(\bfX)
\end{equation}
where $\bfu$ is the displacement of the mean surface and $\bar{\bfu}$ is a director on the mean surface.
The deformation gradient on $\Omega_t$,
\begin{equation}
\bfF := \Gr\,\bfx =\left(\nabla\otimes\bfx \right)^{\rm T} ,
\end{equation}
can then be computed as
\begin{align}
\bfF  = {}&\bfI + \Gr_\Gamma \, \bfu \cdot\Gr\, \bfX_{\Gamma}\nonumber\\ & +D(\bfX)\, \Gr_\Gamma\, \bar{\bfu}\cdot\Gr\bfX_{\Gamma}+\bar{\bfu}\otimes\nabla D(\bfX)\label{Fdef}
\end{align}
where we used the notation
\[
\Gr\, \bff(\bfX_\Gamma(\bfX)) = \Gr_{\Gamma}\, \bff\cdot\Gr\,\bfX_{\Gamma}
\]
and the definitions
\[
\Gr\, \bff := \left(\nabla\otimes\bff\right)^{\rm T} = \frac{\partial \bff}{\partial \bfX},\]\[ \Gr_{\Gamma}\, \bff := \left(\nabla_{\Gamma}\otimes\bff\right)^{\rm T} . 
\]
Note in particular that
\begin{equation}\label{rightproj}
\Gr_{\Gamma}\, \bff = \Gr\, \bff\cdot\bfT_\Gamma .
\end{equation}
From (\ref{closestP}) it follows that
\begin{equation}\label{GradX}
\Gr\, \bfX_{\Gamma}\vert_\Gamma = \bfT_{\Gamma} -D(\bfX)\,\bfkappa(\bfX),
\end{equation}
where
$\bfkappa := \nabla\otimes \bfN$
is the (symmetric) curvature tensor. Using (\ref{Normal}), (\ref{rightproj}), and (\ref{GradX}) in (\ref{Fdef}) we find that
\begin{equation}\label{regularity}
\bfF\vert_{\Gamma}  
=  \bfI + \Gr_\Gamma \,\bfu
+ \bar{\bfu}\otimes\bfN
.
\end{equation}
\subsection{The membrane shell equations}

Consider a potential energy functional given by
\[
\Pi(\bfu,\bar{\bfu}) := \Psi(\bfu,\bar{\bfu})-\Pi^\text{ext}(\bfu,\bar{\bfu})
\]
where $\Psi$ is the strain energy functional and $\Pi^\text{ext}$ is the potential of external loads. We will assume conservative loading so that $\Pi^\text{ext}(\bfu,\bar{\bfu})  = l_\Gamma(\bfu)$ is linear.
Under the assumption of small thickness, we have
\[
\Psi(\bfu,\bar{\bfu}) = \int_{\Omega_t}\Psi_{\bfX}(\bfF)\, d\Omega_t \approx \int_{\Gamma}t\, \Psi_{\bfX}(\bfF) \, d\Gamma =: \Psi_{\Gamma}(\bfu,\bar{\bfu}) ,
\]
where $\Psi_{\bfX}$ is the strain energy per unit volume.
The solution to the nonlinear elastic membrane problem can then be found by seeking stationary solutions to the potential energy functional
\[
\Pi_{\Gamma} := \Psi_{\Gamma}(\bfu,\bar{\bfu}) - l_\Gamma(\bfu) .
\]
Minimizing the potential energy leads to the variational problem of finding $\bar{\bfu}\in [L_2(\Gamma)]^3$, cf. (\ref{regularity}), and $\bfu \in V$, where $V$
is an appropriate Hilbert space which we specify below, such that
\begin{equation}\label{membraneshelleq}
\int_{\Gamma}t \bfP(\bfF) : \left(\nabla_\Gamma\otimes \bfv\right)^{\rm T}\, d\Gamma = l_\Gamma(\bfv) \quad \forall \bfv \in V ,
\end{equation}
and
\begin{equation}\label{Psweak}
\int_{\Gamma}t \bfP(\bfF) : \left(\bar{\bfv}\otimes \bfN\right)\, d\Gamma = 0 \quad\forall \bar{\bfv}\in \left[L_2(\Gamma)\right]^3.
\end{equation}
In (\ref{membraneshelleq}) and (\ref{Psweak}), the first Piola-Kirchhoff stress tensor
\[
 \bfP(\bfF) := \frac{\partial\Psi_{\bfX}}{\partial\bfF}
\]
was introduced in standard fashion. For future use, we shall also define the elastic tangent stiffness as
\[
 \bfL(\bfF) := \frac{\partial^2\Psi_{\bfX}}{\partial\bfF\otimes\partial\bfF}.
\]

The weak form obtained by varying the director field (\ref{Psweak}) can be written on strong form as
 \begin{equation}
\bfP(\bfF)\cdot\bfN = {\bf 0}\quad\text{on $\Gamma$},
\end{equation}
expressing the condition of plane stress in the deformed configuration. We now define the surface deformation gradient
\begin{equation}
\bfF_\Gamma :=  \bfI +\left(\nabla_{\Gamma}\otimes\bfu\right)^{\rm T}
\end{equation}
and the implicit plane stress constitutive relation\footnote{Henceforth, we use the notation $\{\bullet\}$ to indicate an implicit function.}
\begin{equation}\label{Psc}
\bfP_{\Gamma}\{\bfF_\Gamma\} := 
\bfP (\bfF_\Gamma + \bar{\bfu}\otimes\bfN).
\end{equation}
where $\bar{\bfu}$ is chosen so that
\[
\bfP(\bfF_\Gamma +\bar{\bfu}\otimes\bfN)\cdot\bfN = \bfzero.
\]

Using the implicit (plane-stress) function in (\ref{Psc}), we may then state the large deformation membrane shell problem as follows: Find $\bfu\in V$ such that
\begin{equation}
\int_{\Gamma}t \bfP_\Gamma\{\bfF_\Gamma\} : \left(\nabla_{\Gamma}\otimes\bfv\right)^{\rm T}\, d\Gamma = l_\Gamma(\bfv) \quad\forall {\bfv}\in V .
\end{equation}

\subsection{On the plane stress formulation}

From the plane stress condition above, it can be concluded that no traction acts on the surface with normal $\bfN$. Hence, $\bfP$ is an in-plane tensor with respect to its second leg\footnote{Note that the first Piola Kirchhoff stress tensor is a two-point tensor with ``legs'' in both current and reference configurations.} living in the reference configuration. In the current configuration, on the other hand, $\bfP_{\Gamma}$ is an in-plane tensor with respect to its first leg in the current configuration in the sense that
\begin{equation}\label{currps}
\bfn\cdot\bfP(\bfF_\Gamma +\bar{\bfu}\otimes\bfN)= \bfzero,
\end{equation}
where the current normal $\bfn$ can be obtained from the contra-variant transformation
\[
\bfn=\frac{\bfF^{-{\rm T}}\cdot \bfN}{|\bfF^{-{\rm T}}\cdot \bfN|}.
\]
The proof of (\ref{currps}) follows from the fact that the pull-back of the first Piola-Kirchhoff stress to the reference configuration, the second Piola-Kirchhoff stress tensor, is symmetric, which holds for objective, frame-invariant, models of elasticity, cf. \cite{BoWo08}.

Finally, we note that the present formulation of plane stress follows directly from the choice of the kinematical description in (\ref{kinematics}). This can be seen as an alternative to stating the condition in current or material configuration as presented in e.g. \cite{Hughes81} or \cite{Betsch96}, respectively. It should be noted that a variation
\[
    {\rm d}\bfF={\rm d}\bar{\bfu}\otimes\bfN
\]
can never represent a rigid body rotation. Hence, the proposed method is in that sense also adapting the strains to ensure plane stress.

\section{The finite element method\label{FEM}}

\subsection{Parametrization}

In this Section, we explicitly state the matrix representation of the various tensorial quantities needed for the FE implementation.

We assume that we have a shape regular triangulation $\mathcal{T}_h$ of our undeformed midsurface {$\Gamma$}, resulting in a discrete
surface $\Gamma_h$.
For the parametrization of $\Gamma_h$ we wish to define a map $\bfG : (\xi,\eta)\rightarrow \bfX$ from a reference triangle $\hat T$ defined in a local coordinate system $(\xi, \eta)$ to any given triangle $T$ on $\Gamma_h$. Thus the coordinates of the discrete surface are functions of the reference coordinates inside each element,
$\bfX_{\Gamma}^h = \bfX_{\Gamma}^h(\xi, \eta)$.
For any given parametrization, we can extend it to $\Omega_t$ by
defining
\[
\bfX^h(\xi,\eta,\zeta) = \bfX_{\Gamma}^h(\xi, \eta)+\zeta\,\bfN^h(\xi,\eta)
\]
where $-t/2 \leq \zeta \leq t/2$ and $\bfN^h$ is the normal to $\Gamma_h$, found by taking the cross product of $\frac{\partial \bfX_{\Gamma}^h}{\partial \xi}$ and $\frac{\partial \bfX_{\Gamma}^h}{\partial \eta}$.

We consider in particular a finite element parametrization of $\Gamma_h$ as
\begin{equation}\label{surfpara}
\bfX_{\Gamma}^h(\xi,\eta)= \sum_i\bfX_i\psi_i(\xi,\eta)
\end{equation}
where $\bfX_i$ are the physical location of the (geometry representing) nodes on the initial midsurface and $\psi_i(\xi,\eta)$ are finite element shape functions of a certain degree on the reference element.

For the approximation of the displacement, we use a constant extension,
\begin{equation}
\bfu \approx \bfu^h = \sum_i\bfu_i\varphi_i(\xi,\eta)
\end{equation}
where $\bfu_i$ are the nodal displacements, and $\varphi_i$ are shape functions, not necessarily of the same degree as the $\psi_i$.
We employ the usual finite element approximation of the physical derivatives of the chosen basis $\{\varphi_i\}$ on the surface, at $(\xi,\eta)$, in matrix representation,\footnote{For illustration, we here give the explicit component forms in the orhogonal $XYZ$ and $\xi\eta\zeta$ systems, respectively.} as
\[
\left[\begin{array}{>{\displaystyle}c}
\frac{\partial \varphi_j}{\partial X}\\[2mm]
\frac{\partial \varphi_j}{\partial Y}\\[2mm]
\frac{\partial \varphi_j}{\partial Z}\end{array}\right] = \bfJ^{-1}(\xi,\eta,0) \left[\begin{array}{>{\displaystyle}c}
\frac{\partial \varphi_j}{\partial \xi}\\[2mm]
\frac{\partial \varphi_j}{\partial \eta}\\[2mm]
\frac{\partial \varphi_j}{\partial \zeta}\end{array}\right]_{\zeta=0} =: \bfJ^{-1}(\xi,\eta,0)\nabla_{\bfxi}\varphi_j\vert_{\zeta=0},
\]
where $\bfJ(\xi,\eta,\zeta) := \nabla_{\bfxi}\otimes\bfX^h$.
%
This gives, at $\zeta=0$,
\begin{equation}\label{jacobiinv}
\left[\begin{array}{>{\displaystyle}c}
\frac{\partial \varphi_i}{\partial X}\\[2mm]
\frac{\partial \varphi_i}{\partial Y}\\[2mm]
\frac{\partial \varphi_i}{\partial Z}\end{array}\right] = \bfJ^{-1}(\xi,\eta,0) \left[\begin{array}{>{\displaystyle}c}
\frac{\partial \varphi_i}{\partial \xi}\\[2mm]
\frac{\partial \varphi_i}{\partial \eta}\\[2mm]
0\end{array}\right] .
\end{equation}
With the approximate normals we explicitly obtain
\[
\left.\frac{\partial \bfX^h}{\partial \zeta}\right|_{\zeta=0}=\bfN^h ,
\]
so
\[
\bfJ(\xi,\eta,0) := \left[\begin{array}{>{\displaystyle}c>{\displaystyle}c>{\displaystyle}c}
\frac{\partial X^h}{\partial \xi} & \frac{\partial Y^h}{\partial \xi} & \frac{\partial Z^h}{\partial \xi}\\[3mm]
\frac{\partial X^h}{\partial \eta} & \frac{\partial Y^h}{\partial \eta} & \frac{\partial Z^h}{\partial \eta}\\[3mm]
N^h_y & N^h_y & N^h_z\end{array}\right] .
\]
Explicitly we can then write $\nabla \otimes \bfu^h=\sum \nabla\varphi_i\otimes\bfu_i$, and introducting
$\bfT_{\Gamma_h}:= \bfI -\bfN^h\otimes\bfN^h$ we have
$\nabla_{\Gamma_h}\otimes\bfu^h :=\bfT_{\Gamma_h}\cdot\nabla \otimes \bfu^h$.

\subsection{Finite element formulation}

We can now introduce finite element spaces constructed from the basis previously discussed by defining
\begin{equation}\label{spacevA}
W^h_k := \{ v: {v\vert_T \circ\bfG\in P^k(\hat T),\; \forall T\in\mathcal{T}_h};\; v\; \in C^0(\Gamma_h)\},
\end{equation}
and the finite element method reads: Find $\bfu^h\in V^h := [W^h_k]^3$ such that
\begin{equation}\label{membraneFEM}
a_{\Gamma_h}(\bfu^h,\bfv) = l_{\Gamma_h}(\bfv) ,\quad \forall \bfv \in V^h,
\end{equation}
where
\begin{equation}
a_{\Gamma_h}(\bfu,\bfv)   = \int_{\Gamma_h}t \bfP_{\Gamma_h}\{\bfF_{\Gamma_h}\} : \left(\nabla_{\Gamma_h}\otimes\bfv\right)^{\rm T}\, d\Gamma_h .
\end{equation}
Here, the discrete surface deformation gradient is defined as
\begin{equation}
\bfF_{\Gamma_h} :=\bfI +\left(\nabla_{\Gamma_h}\otimes\bfu^h\right)^{\rm T}
\end{equation}
and the relevant implicit plane stress function is
\begin{equation}\label{PShfunc}
\bfP_{\Gamma_h}\{\bfF_{\Gamma_h}\} := \bfP (\bfF_{\Gamma_h} + \bar{\bfu}\otimes\bfN^h)
\end{equation}
where $\bar{\bfu}\in \Bbb{R}^3$ now is solved for locally so that
\begin{equation}\label{PShprob}
\bfP(\bfF_{\Gamma_h} +\bar{\bfu}\otimes\bfN^h)\cdot\bfN^h = \bfzero.
\end{equation}

\subsection{Plane stress iterations}

In order to solve the finite element problem, we adopt Newton iterations on two levels, related to the solution of (\ref{membraneFEM}) and (\ref{PShfunc}), respectively. To this end, we shall now present the algorithmic formulation and the relevant linearizations.

For given $\bfF_{\Gamma_h}$, the local solution $\bar{\bfu}$ to the plane stress problem (\ref{PShprob}) is solved for iteratively as follows: For previous iteration $\bar{\bfu}^{(k)}$, the new iterative solution is defined as $\bar{\bfu}^{(k+1)}=\bar{\bfu}^{(k)}+\Delta\bar{\bfu}$, where the update $\Delta \bar{\bfu}\in\Bbb{R}^3$ is computed from the linear problem
\begin{equation}\label{PSNewt}
    \bfL_{NN}(\bfF_{\Gamma_h} +\bar{\bfu}^{(k)}\otimes\bfN^h)\cdot\Delta\bar{\bfu} = -\bfP(\bfF_{\Gamma_h} +\bar{\bfu}^{(k)}\otimes\bfN^h)\cdot\bfN^h.
\end{equation}
Here, using the notation that $\bfP\cdot\bfN_h=\left[\bfI\otimes\bfN_h\right]:\bfP$, the Jacobian can be expressed as
\begin{equation}\label{LNN}
    \bfL_{NN}(\bfF)= \left[\bfI\otimes\bfN_h\right]:\bfL \cdot \bfN_h,
\end{equation}
where we recall the continuum tangent stiffness $\bfL(\bfF)$. It is easy to show that if $\bfL_{NN}$ is a symmetric tensor\footnote{Note that $\bfL_{NN}$ is of second order.} whenever $\bfL$ possesses major symmetry. Furthermore, $\bfL_{NN}$ will be positive definite, and thus invertible,  as long as the only singular parts of $\bfL$ pertain to rigid body rotations.

The global membrane problem (\ref{membraneFEM}) is solved by finding updates $\Delta \bfu^h \in V^h$, for each previous iteration $\bfu^{h(k)}$, such that
\begin{equation}\label{membraneNewton}
a'_{\Gamma_h}(\bfu^{h(k)},\bfv,\Delta \bfu_h) = l_{\Gamma_h}(\bfv) - a_{\Gamma_h}(\bfu^{h(k)},\bfv),\quad \forall \bfv \in V^h,
\end{equation}
resulting in the iterative solutions $\bfu^{h(k+1)}=\bfu^{h(k)}+\Delta \bfu^h \rightarrow \bfu^h$ with $k$. In (\ref{membraneNewton}), we introduced the tangent form, being the directional derivative of $a_{\Gamma_h}(\bullet,\bullet)$,
\begin{equation}\label{aprime}
a'_{\Gamma_h}(\bfu,\bfv,\bfw) := \int_{\Gamma_h}t \left(\nabla_{\Gamma_h}\otimes\bfv\right)^{\rm T} : \bfL_{\Gamma_h}: \left(\nabla_{\Gamma_h}\otimes\bfw\right)^{\rm T}\, d\Gamma_h,
\end{equation}
where we introduce the plane stress tangent stiffness
\begin{equation}\label{consistentlin}
\bfL_{\Gamma_h} :=  \frac{{\rm d}\bfP_{\Gamma_h}}{{\rm d}\bfF_{\Gamma_h}},
\end{equation}
which is the consistent linearization of the plane stress function defined in (\ref{PShfunc}).

In order to construct the consistent linearization in (\ref{consistentlin}), we study the linearization of (\ref{PShprob}) around the converged solution $\bar{\bfu}$ for given $\bfF_{\Gamma_h}$ as follows:
\begin{equation}
    \bfL_{NN}(\bfF_{\Gamma_h} +\bar{\bfu}\otimes\bfN^h)\cdot{\rm d}\bar{\bfu} +\bfL_N(\bfF_{\Gamma_h} +\bar{\bfu}\otimes\bfN^h):{\rm d}\bfF_{\Gamma_h} = \bfzero.
\end{equation}
We recall the linearization of $\bfP\cdot\bfN_h$ w.r.t. $\bar{\bfu}$ from (\ref{LNN}) and introduce the appropriate linearization w.r.t. $\bfF_{\Gamma_h}$ as
\begin{equation}
    \bfL_{N \bullet}:= \frac{\partial \left[\bfP\cdot\bfN_h\right]}{\partial \bfF_{\Gamma}}=\left[\bfI\otimes\bfN_h\right]:\bfL.
\end{equation}
In the same manner, we identify the linearization of $\bfP$ w.r.t. the director as
\begin{equation}
    \bfL_{\bullet N}:= \frac{\partial \bfP}{\partial \bar{\bfu}}=\bfL\cdot \bfN_h.
\end{equation}
Finally, we may evaluate the sensitivity of $\bar{\bfu}$ w.r.t. $\bfF_{\Gamma_h}$ and (formally) write the consistent linearization of $\bfP_{\Gamma_h}$ as
\begin{equation}\label{statcond}
\bfL_{\Gamma_h} =  \bfL-\bfL_{\bullet N}\cdot\bfL_{N N}^{-1}\cdot\bfL_{N \bullet}.
\end{equation}

\begin{rem}
For a tangent stiffness $\bfL$ satisfying major symmetry, we can show the symmetry
\begin{equation}\label{consistentlin}
    \bfv\cdot \bfL_{N \bullet}:\bfeps \equiv \bfeps:\bfL_{\bullet N}\cdot\bfv \quad \forall (\bfv,\bfeps)\in\Bbb{R}^3\times\Bbb{R}^{(3\times3)},
\end{equation}
i.e., the global tangent stiffness $\bfL_{\Gamma_h}$ is symmetric.
\end{rem}
\section{Constitutive modeling}\label{constitutive}

\subsection{Special case: Linear elasticity}
Linear elasticity can be obtained in the present formulation by simply setting
\begin{equation}\label{linelast}
    \bfP(\bfF)=\bfL_{\rm Hooke}:\left[\bfF-\bfI \right],
\end{equation}
where $\bfL_{\rm Hooke}$ is the constant fourth order constant Hooke tensor related to linear elasticity, satisfying major as well as minor symmetry, cf. the discussion in \cite{LarssonRunesson04}. The minor symmetry, imposing a condition of symmetry of the stress as well as an invariance to the skew symmetric part of the displacement gradient, is a direct consequence of a linearization of {\em any} objective finite deformation constitutive model around $\bfF=\bfI$. Hence, the formulation in (\ref{linelast}) is in complete analogy with stating the small strain response
\[
    \bfsig=\bfL_{\rm Hooke}:\bfeps, \quad \bfeps:=\left(\nabla\otimes\bfu\right)^{\rm sym}.
\]

For linear elasticity, it is standard procedure to define the plane stress equations explicitly. Here, we may of course do the same by computing the (constant) stiffness tensor $\bfL_{\Gamma_h}$ off-line. For isotropic elasticity, the classical result of modifying the Lam\'{e} parameters can of course be reproduced with the procedure described in (\ref{statcond}) for the Hooke tensor.

Adopting a linear elastic model on the form (\ref{linelast}) specializes the present framework to the result in \cite{HaLa14}.

\subsection{Mooney--Rivlin}

In our numerical examples, we use a compressible isotropic Mooney--Rivlin model in which we choose
parameters $E$ and $\nu$, and define
$K=E\nu/(1-\nu^2)$, $\mu=E/(2(1+\nu))$, and $\mu_1=\mu_2 = \mu/2$. Then the Mooney--Rivlin strain energy density is given by
\[
\Psi_{\bfX}(\bfF) :=  \frac12\mu_1  \hat{I}_1 + \frac12 \mu_2 \hat{I}_2 + \frac12 K (J-1)^2
\]
where $J := \text{det}\, \bfF$, $\hat{I}_1 := J^{-2/3}I_1$, and $\hat{I}_2 := J^{-4/3}I_2$, with $I_1$ and $I_2$ the first and second invariants of the left Cauchy--Green tensor $\bfB = \bfF\cdot\bfF^{\rm T}$.

\subsection{Remarks on form finding\label{formfind}}

Trying to minimize the area in the current configuration can be viewed as a particular case of potential energy minimization, cf., e.g.,
Bletzinger et al. \cite{BlWuDaCa05};
a boundary driven problem (zero right-hand side) in which we are led to an isotropic Cauchy stress of the form $\bfsig = s\bfI$, where $s$ is a given constant, corresponding to the
first Piola--Kirchhoff stress $\bfP = s\,\text{det}\,\bfF_\Gamma\bfF_\Gamma^{-\rm{T}}$.  In \cite{BlWuDaCa05} it is proposed to stabilize this approach by
adding a term proportional to $\bfF_\Gamma$ so that
\begin{equation}\label{continuation}
\int_{\Gamma}\left(\lambda \text{det}\,\bfF_\Gamma\bfF_\Gamma^{-\rm{T}} +(1-\lambda) \bfF_\Gamma\right) : \left(\nabla_\Gamma\otimes \bfv\right)^{\rm T}\, d\Gamma = 0\quad\forall \bfv\in V,
\end{equation}
where $\lambda\in [0,1)$ is to be chosen.
After convergence, the domain is then successively updated so that $\bfX_\Gamma$ will refer to the current position. We shall focus here on the choice $\lambda = 0$. Then, as pointed out in \cite{BlWuDaCa05}, the algorithm to solve  consists of one linear step followed by an update of the geometry. Since
\[
\bfF_\Gamma :=  \bfI +\left(\nabla_{\Gamma}\otimes\bfu\right)^{\rm T}
\]
we can write (\ref{continuation}) explicitely as finding $\bfu\in V$Êsuch that
\[
\int_{\Gamma}\left(\bfI +\left(\nabla_{\Gamma}\otimes\bfu\right)^{\rm T}\right) : \left(\nabla_\Gamma\otimes \bfv\right)^{\rm T}\, d\Gamma = 0\quad\forall \bfv\in V.
\]
Writing $\bfu=\bfx_\Gamma-\bfX_\Gamma$, we see that
\[
\nabla_{\Gamma}\otimes\bfu = \nabla_{\Gamma}\otimes\bfx_\Gamma -\bfI
\]
and we are in fact solving for $\bfx_\Gamma\in V$ such that
\begin{equation}\label{fixedpoint}
\int_{\Gamma}\left(\nabla_{\Gamma}\otimes\bfx_\Gamma\right)^{\rm T} : \left(\nabla_\Gamma\otimes \bfv\right)^{\rm T}\, d\Gamma = 0\quad\forall \bfv\in V,
\end{equation}
which corresponds to the strong problem
\[
\Delta_\Gamma \bfx_\Gamma = 0,
\]
where $\Delta_\Gamma$ is the Laplace--Beltrami operator on the surface with coordinates $\bfX_\Gamma$. This is a classical formulation of the minimal surface problem, cf., e.g., \cite{Dz91}.

In a fully discrete method for solving (\ref{fixedpoint}) we thus have a sequence of discrete surfaces $\{\Gamma_h^n\}$
for $n=0,1,2,\ldots$, where $\Gamma_h^0$ refers to the original, given, discrete surface. We then have the scheme: given $n$,
find
$\bfx^{h,n+1}_{\Gamma_h^n}\in W^{h,n}_k$ such that
\begin{equation}\label{fixedpoint2}
\int_{\Gamma_h^n}\left(\nabla_{\Gamma_h^n}\otimes\bfx^{h,n+1}_{\Gamma^n_h}\right)^{\rm T} : \left(\nabla_{\Gamma_h^n}\otimes \bfv\right)^{\rm T}\, d\Gamma = 0\quad\forall \bfv\in W^{h,n}_k,
\end{equation}
and let $\Gamma_h^{n+1}$ be the triangulation defined by the nodal positions $\bfx^{h,n+1}_{\Gamma^n_h}$, set $n \leftarrow n+1$ and repeat until convergence. Here, $W^{h,n}_k$ refers to the discrete space on the triangulation of $\Gamma_h^n$. We note that the proposed solution method corresponds to a fixed point iteration scheme. In contrast, the method proposed by Dziuk
\cite{Dz91} is a viscous relaxation method in which a time derivative is artificially added so that the problem
becomes to solve
\[
\frac{\partial\bfx_\Gamma}{\partial t}-\Delta_\Gamma \bfx_\Gamma = 0,
\]
which is done in \cite{Dz91} using a semi--implicit time stepping scheme until a stationary solution is reached. (By semi--implicit we
here mean that the Laplace--Beltrami operator must be established on the known mesh on $\Gamma_h^{n}$ at the beginning of each time step.)

We finally remark that the form finding problem does not pertain to a plane stress problem, i.e., there is no restriction on the traction $\bfP\cdot\bfN$ acting on the mean surface. 

\section{Numerical examples\label{numerical}}

In the numerical examples below, unless otherwise stated, the computations have  been done using a super--parametric approximation with piecewise linears for the displacements
and piecewise quadratics for the geometry representation.

\subsection{Form finding}

\begin{figure}[h]
\begin{center}
\includegraphics[width=4in]{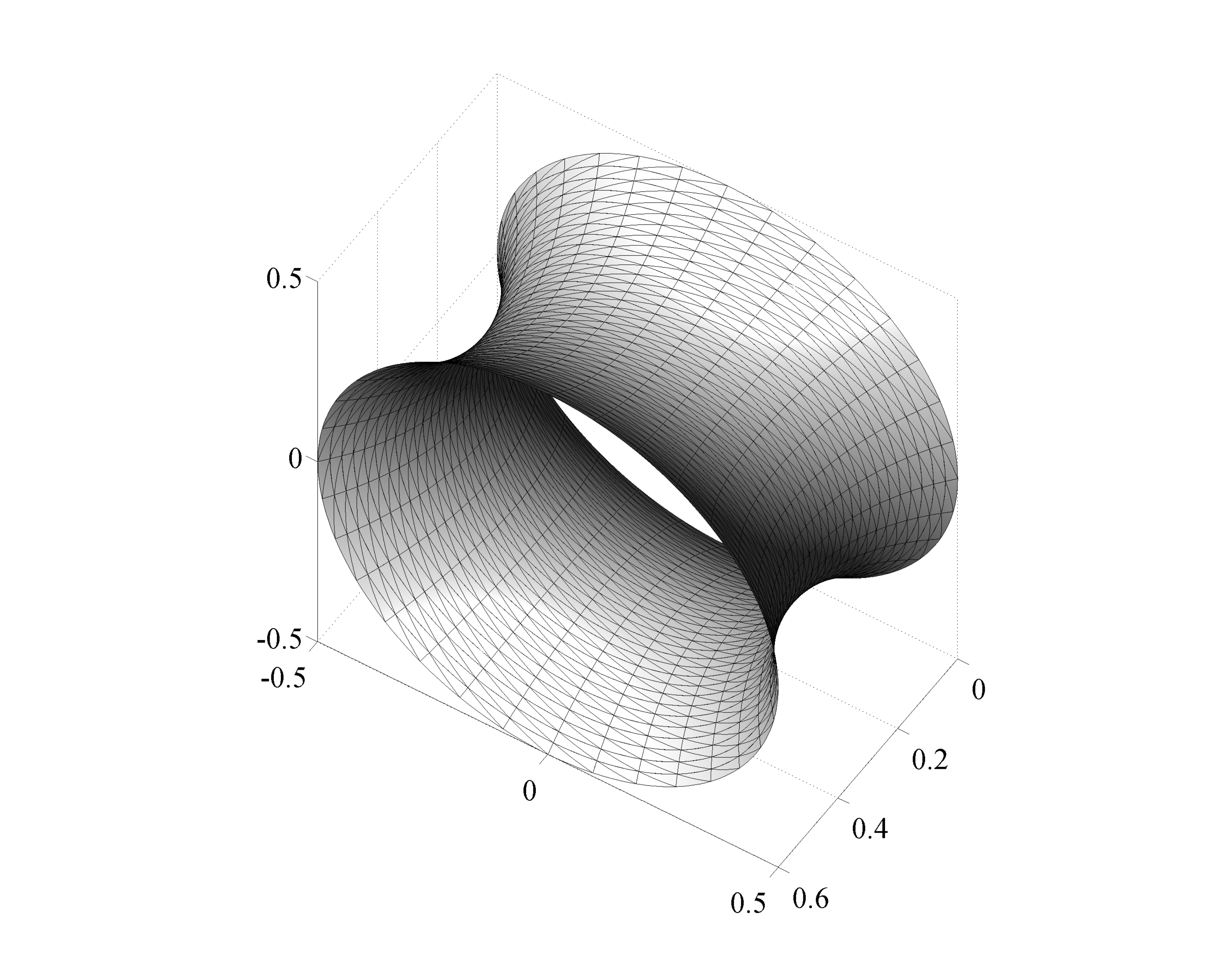}
\end{center}
\caption{Computed catenoid on a particular mesh.}\label{fig:form}
\end{figure}

\begin{figure}[h]
\begin{center}
\includegraphics[width=4in]{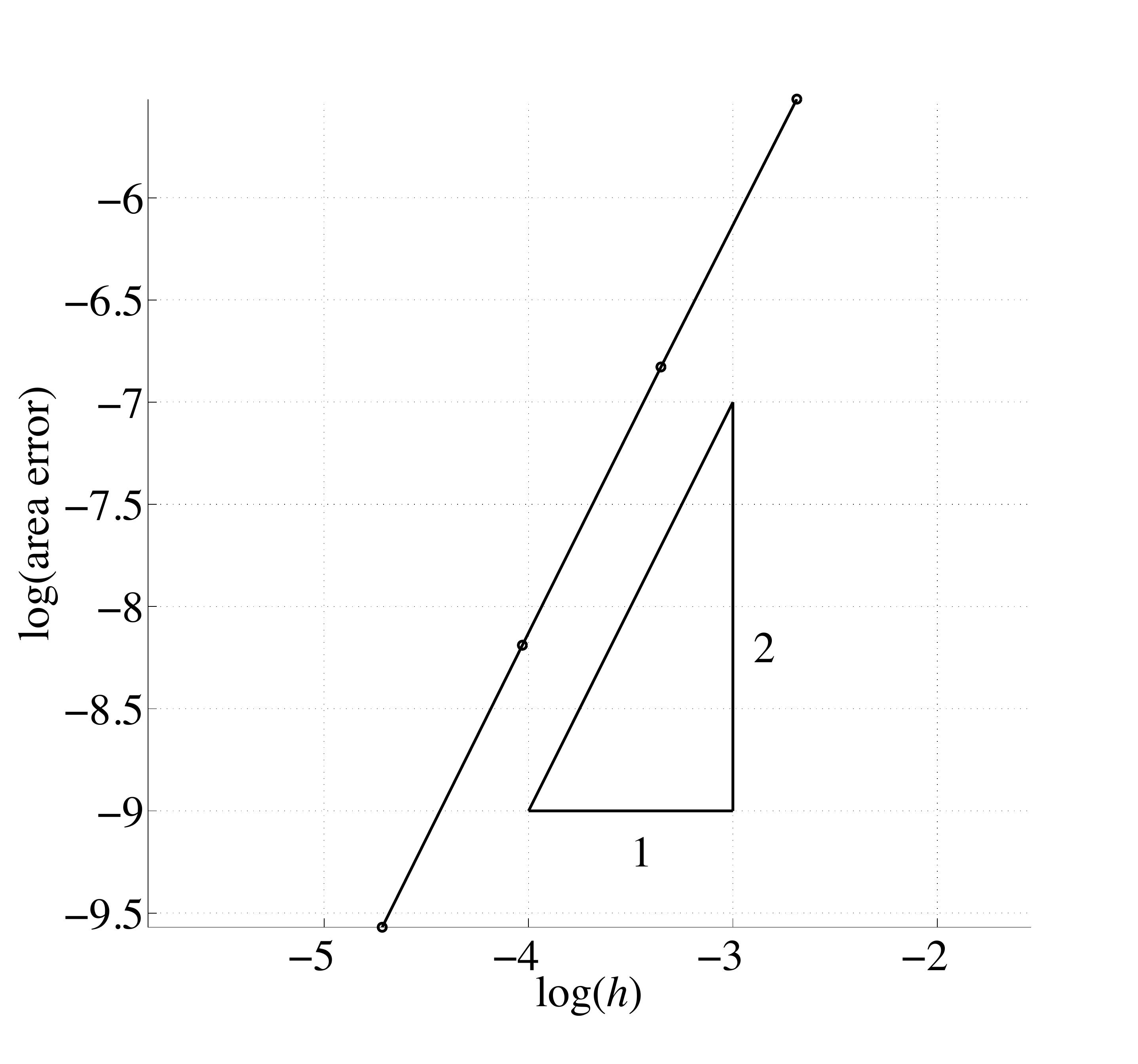}
\end{center}
\caption{Convergence of the area.}\label{fig:formconv}
\end{figure}

We consider a cylinder evolving towards a minimal surface using the algorithm of Section \ref{formfind} with $\lambda = 0$. The
initial radius is 0.5 m, with axis centered at $x=0$, $y=0$, and with height 0.6 m. Both ends of the cylinder are fixed.
The exact solution is then a catenoid whose exact form can be found using Newton's method. In Figure \ref{fig:form} we show
the computed solution on a particular mesh in a sequence of meshes used to check convergence of area, shown if Fig. \ref{fig:formconv}.
We observe second order convergence as expected.

\begin{figure}[h]
\begin{center}
\includegraphics[width=4in]{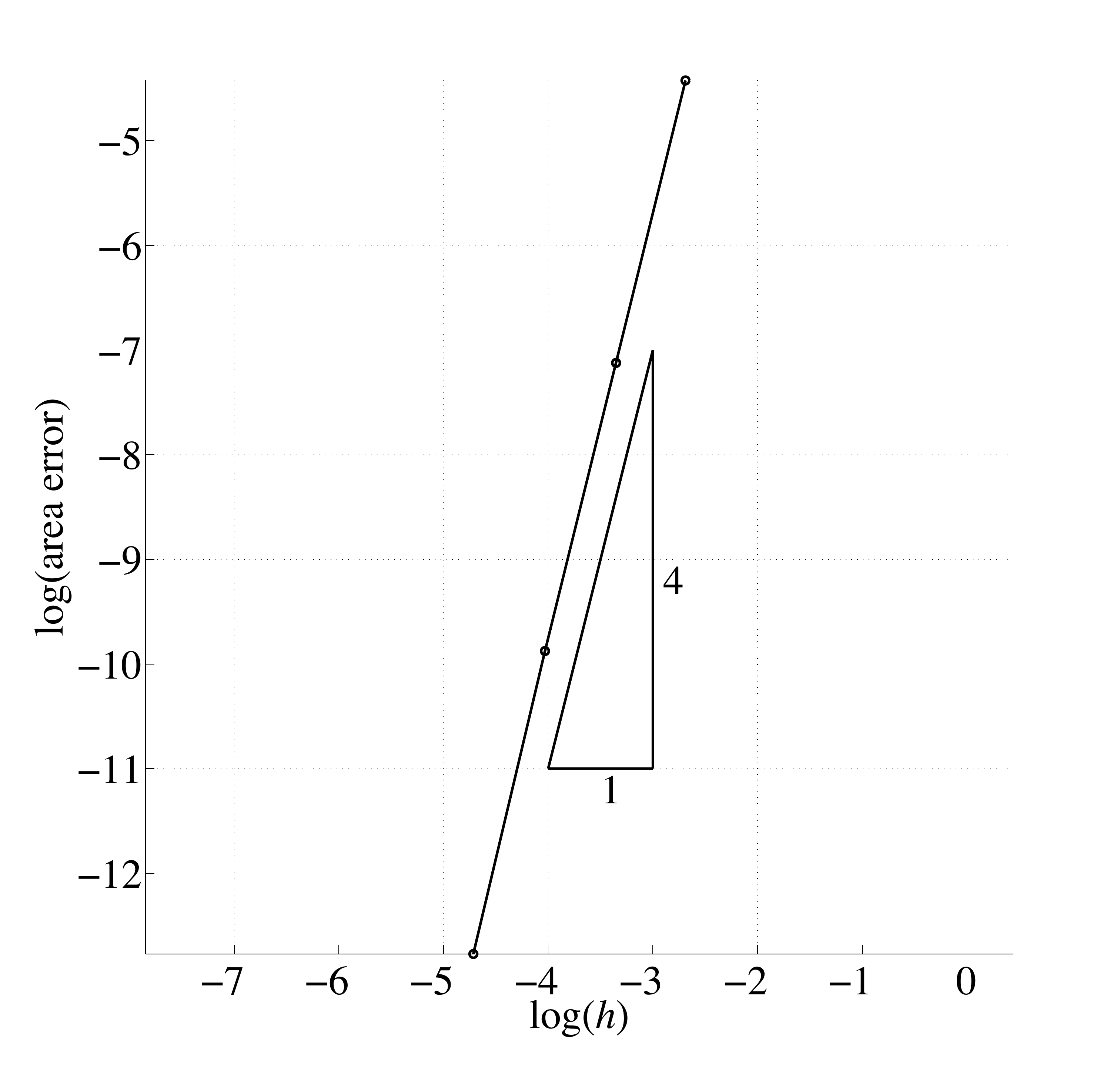}
\end{center}
\caption{Superconvergence of the area, $P^2$ approximation.}\label{fig:super}
\end{figure}

An implementation with isoparametric piecewise quadratic polynomials was also tested and gave superconvergence of order 4 as can be seen in Fig. \ref{fig:super}.
%
This is consistent with the findings of Chien \cite{Ch95}, where it is shown that the 
area of a triangulation (Lagrange--) interpolating a smooth surface is superconvergent for even polynomials but not for odd. 

\subsection{Convergence of $L_2$ norms of the solution}

A cylinder of radius 0.5 m and length 4 m is fixed at both ends and loaded by a conservative force so that
\[
l_{\Gamma_h}(\bfv) := 4000 \int_{\Gamma_h} x\, (4-x)\, \bfN^h\cdot\bfv d\Gamma_h .
\]
The initial thickness was set to $t=1$ cm and the material data were $E=10$ MPa, $\nu=0.5$.

\begin{figure}[h]
\begin{center}
\includegraphics[width=6in]{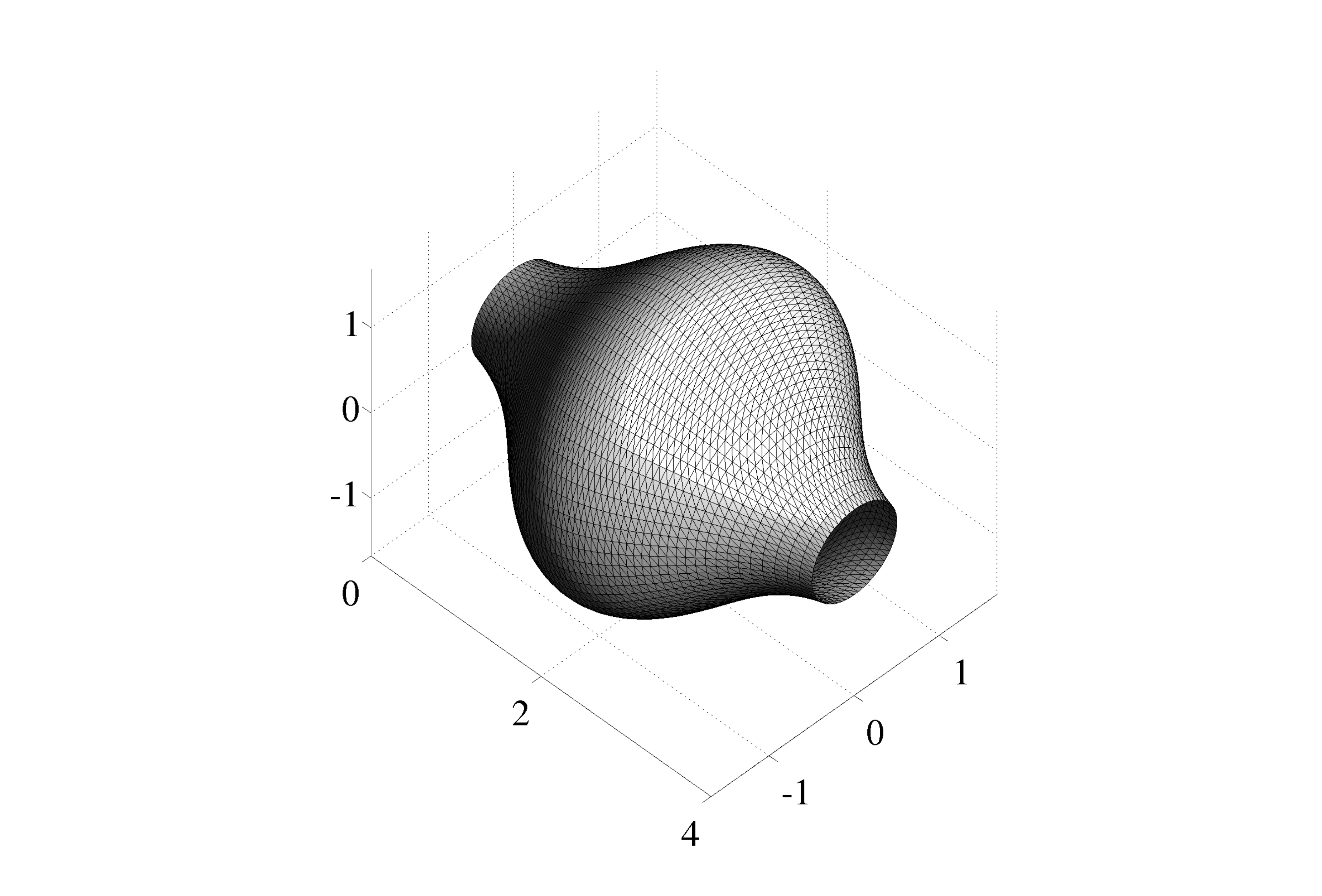}
\end{center}
\caption{Displacements on the finest computational mesh, conservative load.}\label{fig:consdef}
\end{figure}

In Fig. \ref{fig:consdef} we show the deformation on the finest mesh in a sequence for determining convergence.
We check the convergence of norms, so that
\[
e = \| \bfu\|_{L_2(\Gamma)} -\| \bfu^h\|_{L_2(\Gamma_h)} ,
\]
where we replace the exact solution and exact geometry by an ``overkill'' solution (a mesh twice refined from the one in Fig. \ref{fig:consdef}).

\begin{figure}[h]
\begin{center}
\includegraphics[width=4in]{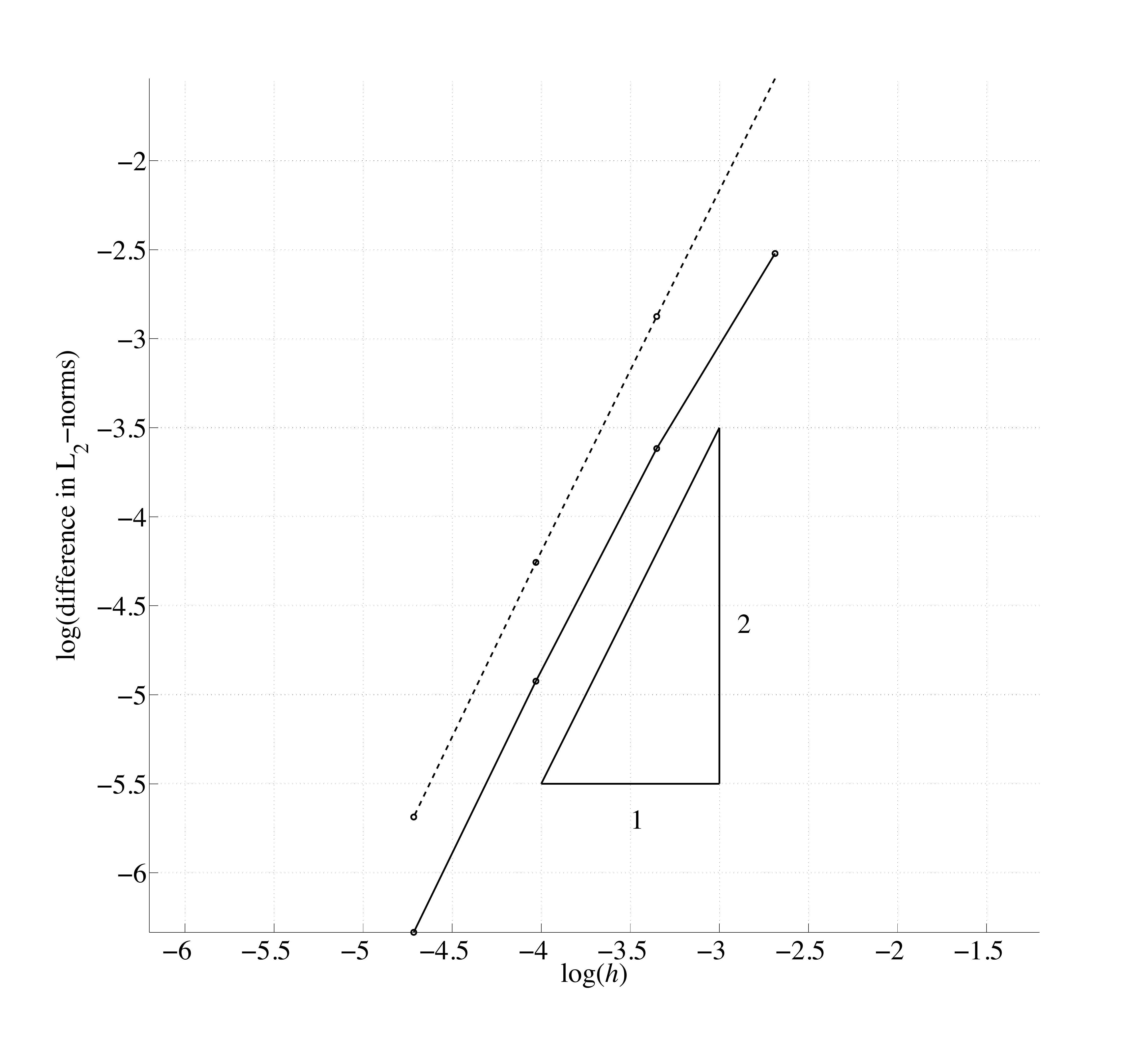}
\end{center}
\caption{Convergence of $L_2$ norms of displacements, conservative load.}\label{fig:consconv}
\end{figure}

In Fig. \ref{fig:consconv} we show the observed convergence of the normal and tangential displacements, both of which show second order convergence.

\subsection{Internal pressure}

In our final example we insider a non-conservative load in the form of an internal pressure. The Newton method then requires linearisation of the
load as discussed by Bonet and Wood \cite[Ch. 8.5.2]{BoWo08}, but in fact we still get convergence without this term, albeit not quadratic convergence.

\begin{figure}[h]
\begin{center}
\includegraphics[width=3in]{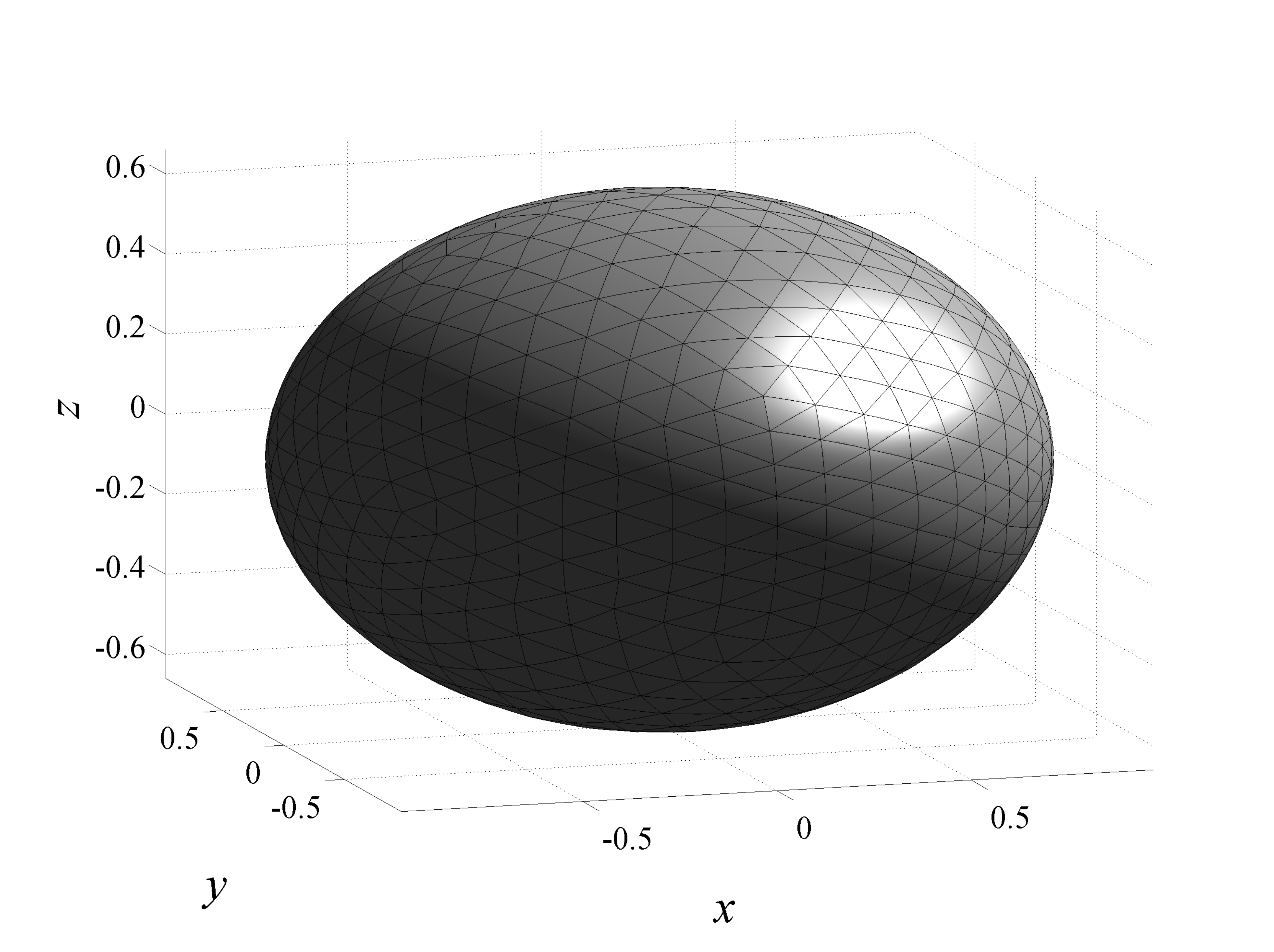}
\end{center}
\caption{Deformation of an oblate spheroid under 1kPa internal pressure.}\label{fig:press1}
\end{figure}

\begin{figure}[h]
\begin{center}
\includegraphics[width=3in]{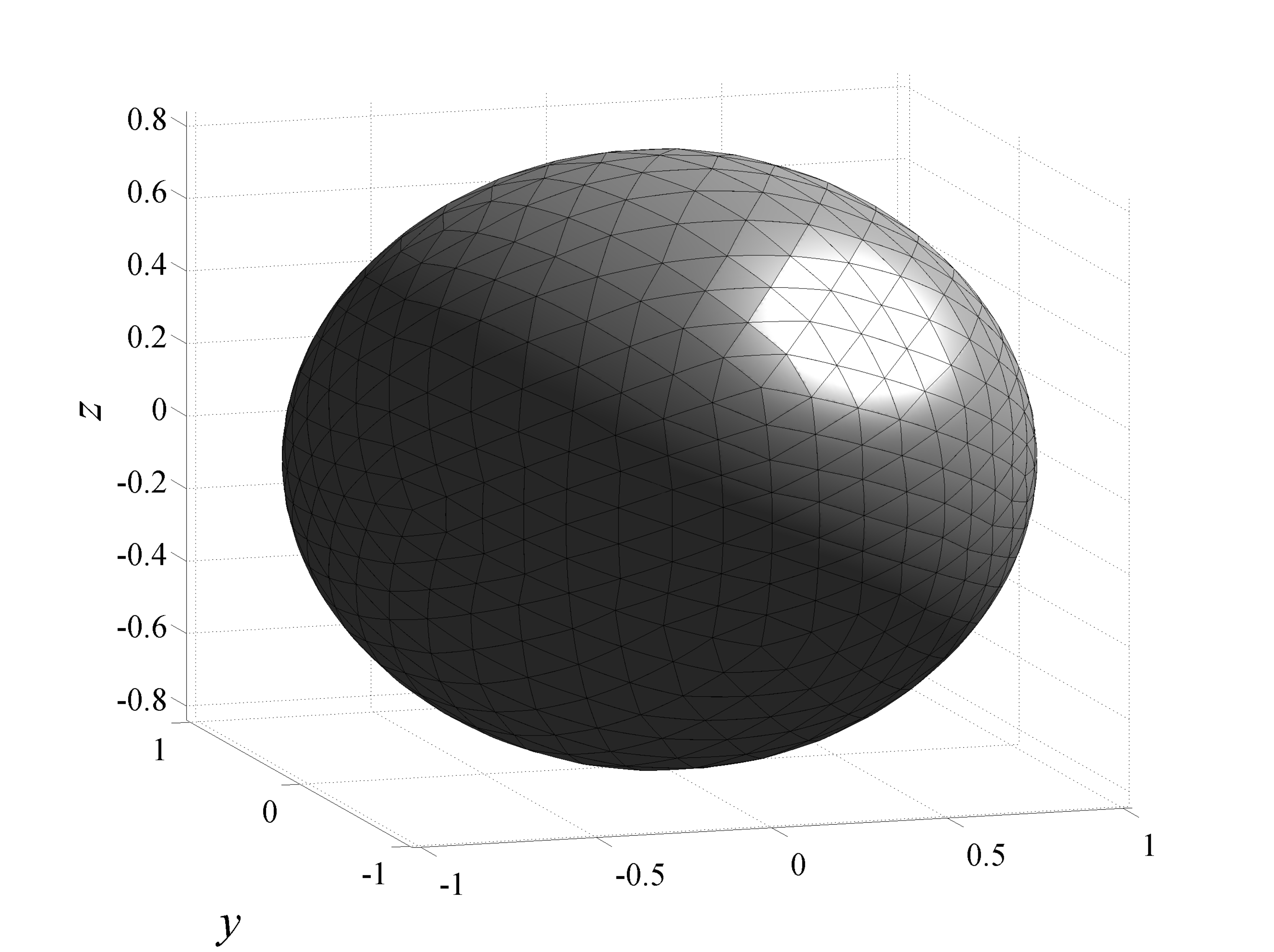}
\end{center}
\caption{Deformation of an oblate spheroid under 3kPa internal pressure.}\label{fig:press2}
\end{figure}

\begin{figure}[h]
\begin{center}
\includegraphics[width=3in]{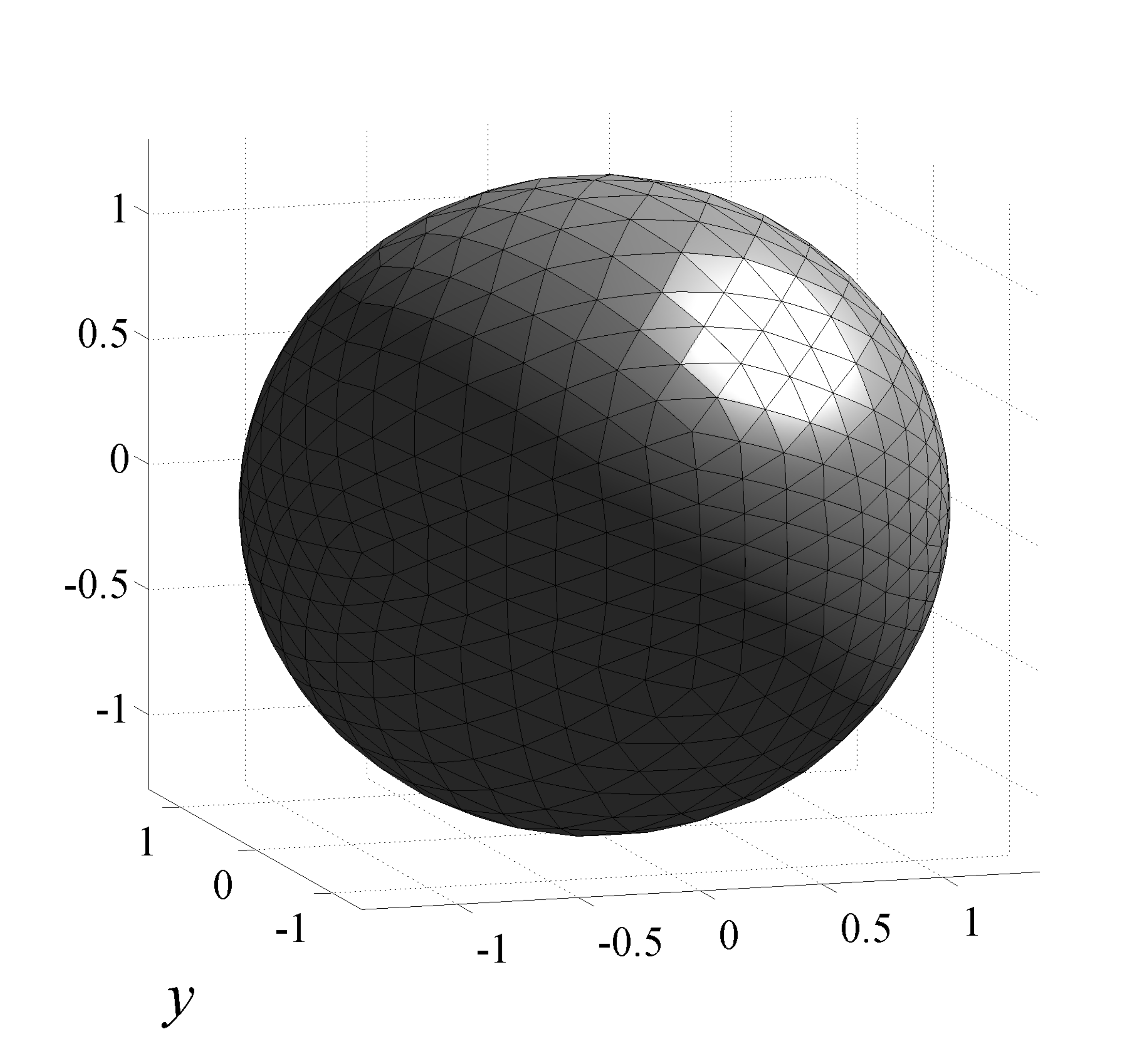}
\end{center}
\caption{Deformation of an oblate spheroid under 4.8kPa internal pressure.}\label{fig:press3}
\end{figure}

\begin{figure}[h]
\begin{center}
\includegraphics[width=4in]{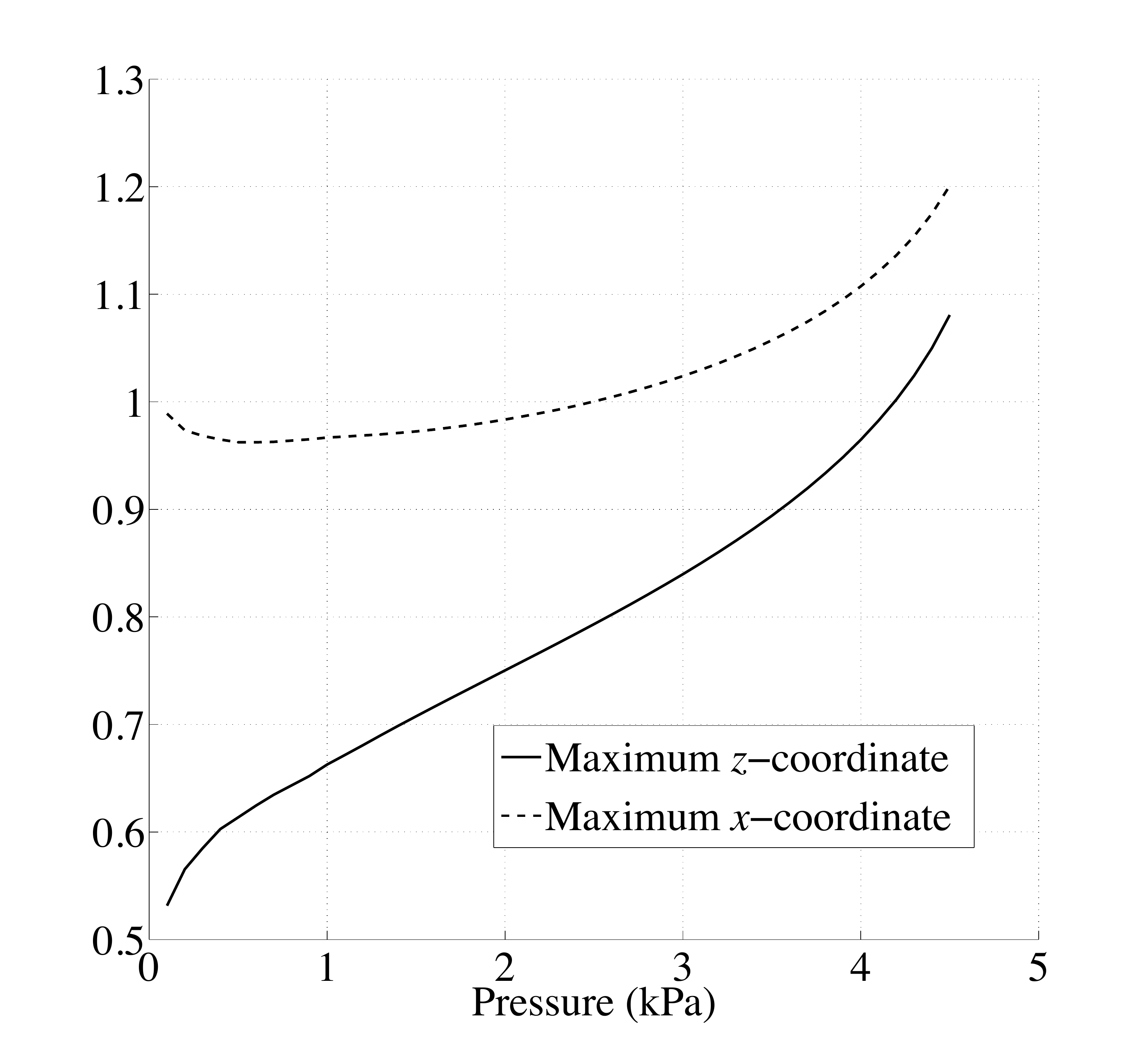}
\end{center}
\caption{Smallest and largest radius of the spheroid under increasing pressure.}\label{fig:radii}
\end{figure}

We show the effect of increasing internal pressure on an oblate spheroid of maximum radius $R_{\max}=1$ m and minimum $R_{\min}=0.5$ m.
The initial thickness is $t =1$ mm and material data are $E=100$ MPa, $\nu=0.5$. In Figs. \ref{fig:press1}--\ref{fig:press3} we show the deformations at 1, 3, and 4.8 kPa,
and in Fig. \ref{fig:radii} we show how the maximum and minimum radii after deformation depend on the pressure. Note that the maximum radius decreases initially.

\section{Concluding remarks}

The finite element method for large deformation membrane elasticity problems
developed herein has the advantage of using Cartesian coordinates as opposed to classical formulations that use co- and contravariant bases.
Our method thus avoids the problem of formulating discrete approximations using these (varying) basis vectors, which may be intricate in a finite element setting if a $C^0$--continuous geometry is used.
Furthermore, our model is continuous---unlike several of the more popular discrete shell models that start out by collapsing 3D continua in a discrete finite element setting, cf. \cite{ChBa03}---which makes modelling more tractable. In addition, a coordinate-invariant formulation of plane stress has also been introduced.

\section*{Acknowledgements}
This research was supported in part by the Swe\-dish Foundation for Strategic Research Grant No.\ AM13-0029 
and the Swe\-dish Research Council Grants No.\ 2011-4992 and No.\ 2013-4708.

\newpage\newpage


\begin{thebibliography}{10}

\bibitem{Betsch96}
P.~Betsch, F.~Gruttmann, and E.~Stein.
\newblock A 4--node finite shell element for the implementation of general
  hyperelastic {3D}--elasticity at finite strains.
\newblock {\em Comput. Methods Appl. Mech. Engrg.}, 130(1--2):57--79, 1996.

\bibitem{BlWuDaCa05}
K.-U. Bletzinger, R.~W{\"u}chner, F.~Daoud, and N.~Camprub{\'{\i}}.
\newblock Computational methods for form finding and optimization of shells and
  membranes.
\newblock {\em Comput. Methods Appl. Mech. Engrg.}, 194(30-33):3438--3452,
  2005.

\bibitem{BoWo08}
J.~Bonet and R.~D. Wood.
\newblock {\em Nonlinear {C}ontinuum {M}echanics for {F}inite {E}lement
  {A}nalysis}.
\newblock Cambridge University Press, Cambridge, second edition, 2008.

\bibitem{ChBa03}
D.~Chapelle and K.~J. Bathe.
\newblock {\em The {F}inite {E}lement {A}nalysis of {S}hells---{F}undamentals}.
\newblock Computational Fluid and Solid Mechanics. Springer-Verlag, Berlin,
  2003.

\bibitem{Ch95}
D.~Chien.
\newblock Numerical evaluation of surface integrals in three dimensions.
\newblock {\em Math. Comp.}, 64(210):727--743, 1995.

\bibitem{Ci00}
P.~G. Ciarlet.
\newblock {\em Mathematical Elasticity. {V}ol. {III}}, volume~29 of {\em
  Studies in Mathematics and its Applications}.
\newblock North-Holland Publishing Co., Amsterdam, 2000.

\bibitem{DeZo95}
M.~C. Delfour and J.-P. Zol{\'e}sio.
\newblock A boundary differential equation for thin shells.
\newblock {\em J. Differential Equations}, 119(2):426--449, 1995.

\bibitem{DeZo97}
M.~C. Delfour and J.-P. Zol{\'e}sio.
\newblock Differential equations for linear shells: comparison between
  intrinsic and classical models.
\newblock In {\em Advances in mathematical sciences: {CRM}'s 25 years
  ({M}ontreal, {PQ}, 1994)}, volume~11 of {\em CRM Proc. Lecture Notes}, pages
  41--124. Amer. Math. Soc., Providence, RI, 1997.

\bibitem{Dz91}
G.~Dziuk.
\newblock An algorithm for evolutionary surfaces.
\newblock {\em Numer. Math.}, 58(6):603--611, 1991.

\bibitem{MoHeYv11}
E.~E.~Monteiro, Q.-C. He, and J.~Yvonnet.
\newblock Hyperelastic large deformations of two-phase composites with
  membrane-type interface.
\newblock {\em Int. J. Eng. Sci.}, 49(9):985--1000, 2011.

\bibitem{LarssonRunesson04}
Larsson F. and Runesson K.
\newblock Modeling and discretization errors in
  hyperelasto--(visco--)plasticity with a view to hierarchical modeling.
\newblock {\em Comput. Methods Appl. Mech. Engrg.}, 193(48--51):5283--5300,
  2004.

\bibitem{HaLa14}
P.~Hansbo and M.~G. Larson.
\newblock Finite element modeling of a linear membrane shell problem using
  tangential differential calculus.
\newblock {\em Comput. Methods Appl. Mech. Engrg.}, 270:1--14, 2014.

\bibitem{HaLaLa14}
P.~Hansbo, M.~G. Larson, and K.~Larsson.
\newblock Variational formulation of curved beams in global coordinates.
\newblock {\em Comput. Mech.}, 53(4):611--623, 2014.

\bibitem{Hughes81}
T.~J.~R. Hughes and W.~K. Liu.
\newblock Nonlinear finite element analysis of shells: {P}art {I}. {T}hree
  dimensional shells.
\newblock {\em Comput. Methods Appl. Mech. Engrg.}, 26(3):331--362, 1981.

\end{thebibliography}
\end{document}